\documentclass[letterpaper, 10 pt, conference]{ieeeconf}  

\IEEEoverridecommandlockouts                              
\newtheorem{theorem}{Theorem}[section]
\newtheorem{lemma}{Lemma}[section]
\newtheorem{proposition}{Proposition}[section]
\usepackage{amsmath, amssymb}
\usepackage{mathtools}
\usepackage{commath}

\DeclareMathOperator*{\argmin}{arg\,min}

\makeatletter
\newcommand{\specialcell}[1]{\ifmeasuring@#1\else\omit$\displaystyle#1$\ignorespaces\fi}
\makeatother

\title{\LARGE \bf
On the Regret of $\mathcal{H}_{\infty}$ Control
}

\author{Aren Karapetyan, Andrea Iannelli, and John Lygeros
\thanks{The authors are with the Department of Information Technology and Electrical Engineering, Automatic Control Laboratory, ETH Zürich, 8092, Zürich, Switzerland. Emails: \{akarapetyan, iannelli, lygeros\}@control.ee.ethz.ch}
\thanks{This work has been supported by the Swiss National Science Foundation under NCCR Automation (grant agreement $51\text{NF}40\_180545$) and under grant no. $200021\_178890$.
}
}

\begin{document}

\maketitle
\thispagestyle{empty}
\pagestyle{empty}

\begin{abstract}
The $\mathcal{H}_{\infty}$ synthesis approach is a cornerstone robust control design technique, but is known to be conservative in some cases. The objective of this paper is to quantify the additional cost the controller incurs planning for the worst-case scenario, by adopting an approach inspired by regret from online learning. We define the \textit{disturbance-reality gap} as the difference between the predicted worst-case disturbance signal and the actual realization. The regret is shown to scale with the norm of this \textit{gap}, which turns out to have a similar structure to that of the certainty equivalent controller with inaccurate predictions, obtained here in terms of the \textit{prediction error} norm.
\end{abstract}

\section{INTRODUCTION}

In this work we focus on the control of linear time-invariant (LTI) systems subject to process noise. The optimal control of such systems for the case of stochastic noise is well studied and is referred to as $\mathcal{H}_2$ control \cite{hassibi1999indefinite}. The controller in this case is optimal for the expected value of the associated cost. When the noise is non-stochastic, a popular approach is to model it as adversarial and having finite energy. In this case, the $\mathcal{H}_{\infty}$ controller  solves the disturbance attenuation problem by ensuring that the system is internally stable and that the infinity norm of the transfer function mapping disturbances to a measurable performance metric is minimized \cite{zhou1998essentials}. For LTI systems with quadratic performance metrics this is equivalent to minimizing the induced $2$-norm of this transfer function matrix \cite{green2012linear}. The closed-loop system is thus robust to any allowable process noise. This robustness proves vital in many applications, especially safety-critical ones. However, as it plans for the worst, $\mathcal{H}_{\infty}$ control can in practice result in a conservative performance and incur a high cost. The aim of of this work is to analyze and quantify the degree of this conservatism using online learning tools \cite{books/daglib/0016248}. By doing so we study how  the concept of robustness from control theory is reflected in a regret formulation.

Recently there has been an increasing interest in quantifying the performance of control algorithms for dynamical systems in terms of regret \cite{hazan2020nonstochastic}. Regret compares the incurred cost of a given (online) algorithm with a clairvoyant one that has full knowledge of the problem. When the latter is restricted to a policy class, its associated regret is referred to as \textit{policy regret}. This is contrasted to \textit{dynamic regret} with no restrictions on the optimal controller \cite{goel2020power}. The regret of $\mathcal{H}_2$ controllers has been extensively studied in the literature \cite{dean2018regret, simchowitz2020naive, goel2020power}.   A number of approaches have  been proposed for linear time varying models and general convex costs that achieve sublinear \textit{policy regret} (with respect to the control horizon length) \cite{agarwal2019online, foster2020logarithmic}. Similar regret bounds are also achieved for cases with model uncertainties using adaptive control techniques \cite{hazan2020nonstochastic, dean2018regret}. The effect of predictions on the performance of the controllers in the form of \textit{dynamic regret} has been studied in the context of receding horizon control, both for fixed \cite{yu2020power} and time varying costs \cite{zhang2021regret}. The performance of a robust receding horizon controller for general costs has been characterized in \cite{muthirayan2021robust} in terms of the disturbance gain. Regret-optimal controllers are considered in \cite{goel2020regret,sabag2021regret} with a classical $\mathcal{H}_\infty$ approach and in \cite{martin2022safe,didier2022system} by adopting the framework of system level synthesis.


A number of online learning algorithms have been compared to $\mathcal{H}_{\infty}$ control in numerical experiments to show its relatively poor performance in the absence of online updates \cite{hazan2020nonstochastic, goel2020regret}. By solving the disturbance attenuation problem, the $\mathcal{H}_{\infty}$ controller provides robustness guarantees on the possible effects of the noise on the system. This is contrasted to many recent online learning-inspired algorithms, whose objectives differ from that of the robust controller. 

In this work we seek to quantify the additional cost that the $\mathcal{H}_{\infty}$ controller incurs due to its cautious nature, thus allowing direct analytical comparisons with other algorithms.  The upper bound of its \textit{dynamic regret} is shown to scale with the norm of the difference between the worst-case and the observed disturbance signals, which we define as the \textit{disturbance-reality gap}.  This is achieved using the game-theoretic formulation of the disturbance attenuation problem, which has been extensively studied in the robust control literature  \cite{green2012linear,bacsar2008h}. To the best of our knowledge, this is the first work to derive regret bounds for the $\mathcal{H}_{\infty}$ control. In addition, it aims to promote the application of game-theoretic concepts in the design of novel online algorithms in control. To put the $\mathcal{H}_{\infty}$ control results in perspective, the regret of a certainty equivalent (CE) controller with an erroneous prediction of the noise signal is also derived. It is shown to be proportional to the square of the norm of the \textit{prediction error}, in agreement with results from \cite{zhang2021regret}. Moreover, when the value of this norm is equal to the norm of the \textit{disturbance-reality gap}, the CE controller attains a lower upper bound as compared to the $\mathcal{H}_\infty$ one. This reinforces  the empirically observed good performance of many ``optimistic'' algorithms \cite{zhang2021regret,muthirayan2021online}. We show that, while similar in structure to that of the CE, the regret of the $\mathcal{H}_{\infty}$ controller has additional terms arising due to a mismatch between its stabilizing and the linear quadratic regulator (LQR) optimal state feedback gains. An exemplifying numerical example is finally provided.

\textit{Notation}: For a matrix $A$ the spectral radius and the spectral norm are denoted by $\rho(A)$, and  $\|A\|$, respectively; $\lambda_{max}(A)$ denotes the maximum and $\lambda_{min}(A)$ the minimum eigenvalue of $A$. For a vector $w \in \mathbb{R}^{n  a}$, $w_{[1:a]}:= [w_1^{\top}, ..., w_{a}^{\top}]^{\top}$ and $\|w\|$ denotes its Euclidean norm. The Kronecker product between two matrices is denoted by $\otimes$. 


\section{PRELIMINARIES}

We consider a LTI system 
\begin{equation}
    x_{t+1} = Ax_t+Bu_t+w_t,
    \label{eq:model}
\end{equation}
with initial state $x_0 \in \mathbb{R}^n$  and known matrices $A\in \mathbb{R}^{n \times n}$ and $B \in \mathbb{R}^{n\times m}$. The control input is denoted by $u_t \in \mathbb{R}^m$, and $x_t,w_t\in\mathbb{R}^n$ are the state and disturbance vectors, respectively. The state is assumed to be  fully observed. The control objective is to minimize the total accumulated cost over a horizon  of length $T$
\begin{equation}
    J^T(u,w;x_0) = x_T^\top Q_Tx_T + \sum_{t=0}^{T-1} x_t^\top Qx_t + u_t^\top Ru_t,
    \label{eq:cost_T}
\end{equation}
where $Q,Q_T \in \mathbb{R}^{n \times n}$ and $R \in \mathbb{R}^{m\times m}$ are design matrices, $u:= [u_0^{\top}, ..., u_{T-1}^{\top}]^{\top} \in \mathbb{R}^{mT}$ and $w:= [w_0^{\top}, ..., w_{T-1}^{\top}]^{\top} \in \mathbb{R}^{nT}$. Moreover, it is assumed that $Q,Q_T \succeq 0$, $R \succ 0$, and the pair $(A,Q)$ is detectable, the pair $(A,B)$ is  stabilizable, and  $\|x_0\| \leq X$ for some $X \in \mathbb{R}^+$. We consider finite energy disturbance signals \cite{green2012linear, bacsar2008h} in the space  $\mathcal{L}_{2}(0,T) =\{w:N\rightarrow \mathbb{R}^n \; : \; \|w\|_{[0,T-1]}<\infty\}$ where $\|w\|_{[0,T-1]} =\left(\sum_{k=0}^{T-1}\|w(k)\|^2\right)^{\frac{1}{2}}$ over a finite horizon and $ \mathcal{L}_{2}(0, \infty) = \{w : \|w\|< \infty\}$ over an infinite horizon\footnote{Arguments and subscripts are dropped from this notation when clear from the context}. Disturbance signals considered in this paper will generally take values in the space of $\mathcal{L}_2$ signals with total energy less than or equal to $1$, denoted by $\bar{\mathcal{L}}_2$. This allows us to define the  infinite horizon cost  $J(u,w;x_0) := \lim_{T \rightarrow \infty}J^T(u,w;x_0)$. 

\subsection{The $\mathcal{H}_{\infty}$ Problem}

The robust $\mathcal{H}_{\infty}$ controller  minimizes the induced spectral norm of operator $T_{\mu}$ mapping the disturbance signal $w \in \bar{\mathcal{L}}_2$ to an output signal $z \in {\mathcal{L}}_2$ and internally stabilizing the system for the infinite horizon case \cite{green2012linear}. The optimization problem can  be written as
\begin{equation}
    \adjustlimits\inf_{\mu \in \mathcal{M}}\sup_{w \in \bar{\mathcal{L}}_2}\frac{\|T_{\mu}w\|}{\|w\|} = \adjustlimits\inf_{\mu \in \mathcal{M}} \sup_{\|w\| = 1}\|T_{\mu}w\| : = \gamma^{\star},
\end{equation}
where $\mathcal{M}$ is the set of policies with access to current and past state measurements \cite{bacsar2008h}. With an appropriate definition of the output signal $z$, the above can be written in terms of the total cost\footnote{We abuse the notation slightly, $J(\mu,w;x_0)$ denotes here the cost
associated with the control inputs generated by the policy $\mu$}
\begin{equation}
    \gamma^{\star} = \adjustlimits\inf_{\mu \in \mathcal{M}}\sup_{\|w\| = 1}\left(J(\mu,w;x_0)\right)^{\frac{1}{2}},
    \label{eq:h_inf_problem}
\end{equation}
or with $J$ replaced by $J^T$ for the finite horizon case \cite{bacsar2008h}. A policy $\mu^{\star}$ that attains the minimum in \eqref{eq:h_inf_problem} will be referred to as the $\mathcal{H}_{\infty}$ controller.

\subsection{Regret Definition}
\label{sec:regret_definition}
The regret is a metric designed to measure the performance of online learning algorithms \cite{books/daglib/0016248, hazan2019introduction}. The online convex optimization (OCO) setting considers a decision variable $x_k$, chosen from a given convex set $\mathcal{X}$ at timestep $k$. Cost $c_k$ suffered by the decision maker is then revealed, according to an unknown convex function from within a class $ \mathcal{C}: \mathcal{X} \rightarrow \mathbb{R}$. An algorithm $\mathcal{A}$ maps the available history of cost functions  up to time $k-1$ to a decision variable at time $k$
\begin{equation}
    x_k^{\mathcal{A}}= \mathcal{A}(c_0,\hdots,c_{k-1})\in \mathcal{X}.
\end{equation}
The regret of algorithm $\mathcal{A}$ at timestep $T$ is then defined as
\begin{equation*}
    \mathcal{R}_T^{\mbox{OCO}}(\mathcal{A}) = \sup_{\{c_0,\hdots, c_T\}\subseteq \mathcal{C}}\left(\sum_{k=0}^Tc_k(x_k^{\mathcal{A}})-\min_{x \in \mathcal{X}}\sum_{k=0}^{T}c_k(x_k)\right),
\end{equation*}
which compares the performance of the given algorithm to the one with full knowledge of the unknown. This idea has been extended to dynamical systems subject to process noise \cite{hazan2020nonstochastic}. For non-stochastic, adversarial disturbance signals, the worst-case \textit{policy regret} of an algorithm $\mathcal{A}$ is defined as,
\begin{equation}
    \mathcal{R}^{\Pi}_T(\mathcal{A}) = \sup_{w}\left(\sum_{k=0}^Tg_k(u_k^{\mathcal{A}}, x_k)-\sum_{k=0}^{T}g_k(u_k^{\Pi},x_k^{\Pi})\right),
    \label{eq:regret_dynamic}
\end{equation}
where the system evolves according to a model $x_{k+1} = f_k(x_k,u_k,w_k)$, $g_k$ is the stage cost, the noise is bounded $\|w_k\|\leq W \quad \forall 0\leq k <T$, $u_k^{\mathcal{A}}$ is the control input at timestep $k$ generated by the algorithm $\mathcal{A}$, and $u^{\Pi}$ and $x^{\Pi}$ are the optimal offline control inputs and states, respectively. The optimal offline inputs are the solution of the following optimization problem,
\begin{equation*}
\begin{split}
    u^{\Pi} = &\argmin_{u \in \Pi}\sum_{k=0}^{T}g_k(u_k,x_k)\\
    s.t. \quad &x_{k+1} = f_k(x_k,u_k,w_k),
    \end{split}
\end{equation*}
thus chosen in hindsight given full knowledge of the noise realization and the cost function sequence. Here the optimal inputs are selected from a fixed set of polices $\Pi$. When the inputs are not restricted to any class of policies, \eqref{eq:regret_dynamic} is instead referred to as \textit{dynamic regret}. In this paper we consider only the latter, which for a given noise realisation $w\in \bar{\mathcal{L}_2}$ and for the problem at hand is as follows
\begin{equation}
    \mathcal{R}_T(\mathcal{A},w) = J^T(u^{\mathcal{A}},w;x_0)-J^T(u^{\star},w;x_0).
\end{equation}
Here the optimal offline inputs $u^{\star}$ are the solution of the following problem
\begin{equation*}
\begin{split}
    u^{\star} = &\argmin_{u}J^T(u,w;x_0)\\
    s.t. \quad &x_{t+1} =  Ax_t+Bu_t+w_t.
    \end{split}
\end{equation*}
\addtolength{\textheight}{-0.5cm}   

\section{REGRET BOUNDS FOR $\mathcal{H}_{\infty}$ CONTROL}

\subsection{The Worst-Case Disturbance}
\label{sec:worst_case_disturbance}
In this section, the disturbance signal that attains the highest cost for the $\mathcal{H}_{\infty}$ controller, $\mu^{\star}$, is characterized. For this we adopt the zero-sum game formulation \cite{bacsar2008h, green2012linear} that makes use of the game theoretical toolkit to derive the optimal controller, and is also relevant for getting regret bounds.

For a given $\gamma > 0$, let $J_{\gamma}^T(u,w,x_0)$ be defined as follows
\begin{equation*}
    J^T_{\gamma}(u,w;x_0) = J^T(u,w;x_0) - \gamma^2 \sum_{t=0}^{T-1}\|w_t\|^2, 
\end{equation*}
with the infinite horizon case defined as  $T\rightarrow \infty$. Signals $u$ and $w$ are considered to be two adversarial players in the game trying to respectively minimize or maximize the cost $J_{\gamma}^T(u,w;x_0)$. The  following \textit{min-max} inequality defines the upper and lower values of the game
\begin{equation*}
    \adjustlimits\inf_{\mu \in \mathcal{M}}\sup_{w \in \bar{\mathcal{L}}_2}\{J_{\gamma}^T(\mu,w;x_0)\}^{\frac{1}{2}} \geq \adjustlimits\sup_{w\in \bar{\mathcal{L}}_2}\inf_{\mu \in \mathcal{M}}\{J_{\gamma}^T(\mu,w;x_0)\}^{\frac{1}{2}},
\end{equation*}
where the upper value is, effectively, the soft constrained version of the original problem \eqref{eq:h_inf_problem}. If there exists a policy pair $(u^{SP},w^{SP})$ such that the lower and upper values are equal then it constitutes a saddle-point (SP) solution. The cost $J_{\gamma}^{T, SP}$ that these policies attain is called the value of the game. 
\subsubsection{Finite horizon}
Consider the following condition on $\gamma$,
\begin{equation}
    \Xi = \gamma^2I -M_{t+1}(\gamma) > 0,\qquad \forall t \in [0, T-1],
    \label{eq:gamma_set}
\end{equation}
where $M_t(\gamma)$ is a sequence of matrices generated by the following coupled generalised Riccati equations
\begin{equation}
\begin{aligned}
    \Lambda_t(\gamma) &= I +(BR^{-1}B^T-\gamma^{-2}I)M_{t+1}(\gamma), \\
    M_t(\gamma) &= Q + A^TM_{t+1}(\gamma)\Lambda_t^{-1}(\gamma)A, 
    \label{eq:M_t}
\end{aligned}
\end{equation}
for all $t \in [0, T-1]$ and with $ M_T(\gamma) = Q_T$.  For the case with perfect state measurements, if $\gamma$ satisfies \eqref{eq:gamma_set}, then $J^T_\gamma$ is strictly convex in $u$ and strictly concave in $w$ \cite{bacsar2008h, green2012linear}, and $\Lambda_t(\gamma)$ is nonsingular \cite{bacsar1998dynamic}. In this case the game has a SP solution, and its value is equal to $x_0^TM_0(\gamma)x_0$ \cite{bacsar2008h}. The SP policies can be formulated as a feedback on the current state, and are defined for all $ 0\leq t<T$ by
\begin{align}
u^{\infty}_t &= -R^{-1}B^TM_{t+1}(\gamma)\Lambda_t^{-1}(\gamma)Ax_t := -K_t^{\infty}x_t,
\label{eq:control_strategy}
\\
w^{\infty}_t &= \gamma^{-2}M_{t+1}(\gamma)\Lambda_t^{-1}(\gamma)Ax_t .
\label{eq:disturbance_policy}
\end{align}
As shown in \cite{bacsar1998dynamic}, if both players play the optimal strategies, then the resulting system will evolve according to 
\begin{equation}
    x^{\infty}_{t+1} = \Lambda^{-1}_t(\gamma)Ax^{\infty}_{t}, \quad x_0^{\infty} = x_0 \quad \forall 0\leq t<T.
    \label{eq:feedback_matrix}
\end{equation}
It then follows from the ordered interchangeability property of zero sum games that the disturbance policy \eqref{eq:disturbance_policy} with $x_t$ replaced by $x^{\infty}_t$ from \eqref{eq:feedback_matrix}
\begin{equation}
    w^{\star}_t: = \gamma^{-2}M_{t+1}(\gamma)\Lambda_t^{-1}(\gamma)Ax_t^{\infty} \qquad  \forall 0\leq t<T
\label{eq:disturbance_policy_ol}
\end{equation}
constitutes an open loop strategy in saddle point equilibrium with \eqref{eq:control_strategy}.
Moreover, \cite{limebeer1989discrete, bacsar1998dynamic} show that the feedback policy \eqref{eq:control_strategy} with $\gamma = \underline{\gamma}$, the lowest possible value that satisfies \eqref{eq:gamma_set}, and with initial state $x_0 = 0$, is the $\mathcal{H}_{\infty}$ minimax controller. For this value of $\gamma$ the optimal disturbance attenuation problem \eqref{eq:h_inf_problem} and its soft constrained version coincide, and $\underline{\gamma}= \gamma^{\star}$. For any other fixed $\gamma>\gamma^{\star}$ the solution becomes suboptimal. In \cite{didinsky1992design} the case for non-zero initial states is considered and it is shown that a saddle point solution also exists for the original hard-constrained problem \eqref{eq:h_inf_problem}, given that the energy of the disturbance is at its maximum. Here we only consider the initial states $x_0$ that belong to the following set
\begin{equation}
    X_s = \bigcup_{\gamma>\underline{\gamma}}\{x_0 \in \mathbb{R}^n :  \|w^{\star}(\gamma)\| = 1\},
    \label{eq:x0_set}
\end{equation}
such that the maximum energy is achieved with the policy \eqref{eq:disturbance_policy_ol}. Thus, for all initial states in the set $X_s$, a pure strategy saddle point exists for the $\mathcal{H}_{\infty}$ problem. The optimal strategies are given by \eqref{eq:control_strategy} and \eqref{eq:disturbance_policy_ol} with $\gamma =  \bar{\gamma}$ satisfying $ \|w^{\star}(\bar{\gamma})\| = 1$, i.e. having the maximum allowable energy.

\subsubsection{Infinite horizon}
In the infinite horizon case, the minimal non negative-definite, stationary solution to \eqref{eq:M_t} is considered. In this case, \eqref{eq:M_t} become coupled generalized algebraic Ricatti equations (ARE-s) with solutions $M(\gamma)$ and $\Lambda(\gamma)$. For the disturbance attenuation problem \eqref{eq:h_inf_problem}, a saddle point equilibrium exists also in this case and is given for all $ 0\leq t<T$ by
\begin{align}
u_t^{\infty} &= -R^{-1}B^TM(\bar{\gamma})\Lambda^{-1}(\bar{\gamma})Ax_t := -K^{\infty}x_t,
\label{eq:control_strategy_inf}
\\
w_t^{\infty} &= \bar{\gamma}^{-2}M(\bar{\gamma})\Lambda^{-1}(\bar{\gamma})Ax_t,
\label{eq:disturbance_policy_inf}
\end{align}
where $\bar{\gamma}$ is obtained through a trial-and-error method to satisfy $ \|w^{\star}(\bar{\gamma})\|^2 = 1$, along with certain conditions that allow the minimization problem to be well-posed \cite{ mageirou1976values, willems1971least}. The worst-case open loop disturbance signal can then be defined similar to the finite horizon case
\begin{equation}
    w^{\star}_t: = \bar{\gamma}^{-2}M(\bar{\gamma})\Lambda^{-1}(\bar{\gamma})Ax_t^{\infty} \qquad \forall t\geq 0,
    \label{eq:w_worst_inf}
\end{equation}
where
\begin{equation}
x^{\infty}_{t+1} = \Lambda^{-1}(\bar{\gamma})Ax^{\infty}_{t}, \quad x_0^{\infty} = x_0 \quad \forall t \geq 0.
\end{equation}

\subsection{Regret Analysis}

In this section an upper bound for the regret of the $\mathcal{H}_{\infty}$ problem is obtained. We consider first the infinite horizon case with the minimax controller \eqref{eq:control_strategy_inf}, then the finite one with the controller \eqref{eq:control_strategy}. For a given control input $u$ and disturbance signal $w$, at a generic timestep $i$ (with $0\leq i<T$)  the cost-to-go function is defined as
\begin{equation}
    J_i^T(u^{\infty},w;x_i) = x_T^\top Q_Tx_T + \sum_{t=i}^{T-1} x_t^\top Qx_t + u_t^\top Ru_t,
    \label{eq:cost_to_go}
\end{equation}
and  $J(u,w;x_i) : =  \lim_{T \rightarrow \infty}J_i^T(u,w;x_i)$.\\

\subsubsection{Infinite horizon}
The following result is introduced  to characterize the cost in the infinite horizon case.
\begin{lemma}
For all $F$ with $\rho(F)<1$, $P_i \succ 0 \quad \forall i \geq 0$ and $Q \succ 0$, the iteration $P_{i+1} = F^{\top}P_iF+Q$ converges to a unique value $P$.
\label{the:lyapunov}
\end{lemma}
\textit{Proof}: For any $F$ and $\epsilon_1 > 0$, there exists a matrix norm $\norm[1]{\cdot}$, such that $\norm[1]{F}\leq \rho(F)+\epsilon_1$ \cite{dunford1963spectral}. We define $S_a = \sum_{k=0}^{a}(F^{\top})^kQF^k$. Clearly, as $\rho(F)<1$, for all $\epsilon_2 >0$, we can always find $b>a \in \mathbb{N}$, such that
\begin{equation*}
    \norm[1]{S_a-S_b}\leq \norm[1]{Q}\left(\norm[1]{F}^{2(a+1)}+\dots+\norm[1]{F}^{2b}\right) < \epsilon_2.
\end{equation*}
$S_a$ is therefore a Cauchy sequence in the corresponding Banach space and therefore converges. The iteration can then be written as,
\begin{equation*}
\lim_{i\rightarrow\infty}P_{i+1}= \lim_{i\rightarrow\infty}\left( (F^{\top})^iP_0F^i + \sum_{k=0}^{i}(F^{\top})^kQF^k \right). 
\end{equation*}
The first term on the right hand side converges to $0$, thus the iteration converges to the Lyapunov equation with a unique positive definite solution $P$.\hfill $\square$

Following \cite{goel2020power, yu2020power}, we claim that $J(u^{\infty},w;x_i)$ can be expressed in the form of an extended quadratic function, as formulated in the following lemma.
\begin{lemma}
The infinite horizon cost of the $\mathcal{H}_{\infty}$ controller \eqref{eq:control_strategy_inf}, solving the disturbance attenuation problem \eqref{eq:h_inf_problem}, is given by $J(u^{\infty},w;x_i) = x_i^{\top}P^{\infty}x_i + x_i^{\top}v^{\infty}_i +q^{\infty}_i$, with $P^{\infty} \in \mathbb{R}^{n \times n}, v^{\infty}_i \in \mathbb{R}^n, q^{\infty}_i \in \mathbb{R}\quad \forall i\geq 0$  given in \eqref{eq:cost_to_go_infty}.
\label{the:induction}
\end{lemma}

\textit{Proof}: The finite horizon cost for the controller \eqref{eq:control_strategy_inf} is claimed to be given by $J_i^T(u^{\infty},w;x_i) = x_i^{\top}P_i^{\infty}x_i + x_i^{\top}v^{\infty}_i +q^{\infty}_i$, with some $P^{\infty}_i \in \mathbb{R}^{n \times n}, v^{\infty}_i \in \mathbb{R}^n, q^{\infty}_i \in \mathbb{R}$. Indeed, for $i =T$ this holds trivially, with $P^{\infty}_T = Q_T$ and $v^{\infty}_T,q^{\infty}_T = 0$. Then if the claim holds at $i+1$ the cost-to-go at $i$ satisfies
\begin{align}
\begin{split}
   &J_i^T(u^{\infty},w;x_i) =  x_i^{\top}Qx_i + u^{\infty  \top}_iRu^{\infty}_i \\
   &+ (Ax_i+Bu^{\infty}_i+w_i)^{\top}P^{\infty}_{i+1}(Ax_i+Bu^{\infty}_i+w_i)\\
   &+ (Ax_i+Bu^{\infty }_i+w_i)^{\top}v^{\infty}_{i+1} + q^{\infty}_{i+1}\\
   &= u^{\infty  \top}_i(R+B^{\top}P^{\infty}_{i+1}B)u^{\infty}_i\\
   &+2u^{\infty  \top}_iB^{\top}(P^{\infty}_{i+1}Ax_i +P^{\infty}_{i+1}w_i+\frac{v^{\infty}_{i+1}}{2})\\
   &+x_i^{\top}Qx_i+(Ax_i+w_i)^{\top}P^{\infty}_{i+1}(Ax_i+w_i)\\
   &+(Ax_i+w_i)^{\top}v^{\infty}_{i+1}+q^{\infty}_{i+1}.
    \label{eq:valuefunction}
    \end{split}
\end{align}
Substituting $u^{\infty}_i = -K^{\infty}x_i$, defining $F^{\infty} = A -BK^{\infty}$ and grouping terms leads to
\begin{align}
    \begin{split}
        &J^T(u^{\infty};x_i) = x_i^{\top}\underbrace{(F^{\infty  \top}P^{\infty}_{i+1}F^{\infty} + Q + K^{\infty  \top}RK^{\infty})}_\text{$P_i^{\infty}$}x_i +\\
        & x_i^{\top}\underbrace{(F^{\infty \top}(2P^{\infty}_{i+1}w_i + v^{\infty}_{i+1}))}_\text{$v_i^{\infty}$}+ \underbrace{w_i^{\top}v^{\infty}_{i+1} + w_i^{\top}P^{\infty}_{i+1}w_i+q^{\infty}_{i+1}}_\text{$q_i^{\infty}$}.
    \end{split}
\label{eq:cost_to_go_infty}
\end{align}
The expression for $v_i^{\infty}$ can  be rewritten as
\begin{equation*}
v^{\infty}_i  = 2 \sum_{j=0}^{T-i-1}(F^{\infty  \top})^{j+1}P^{\infty}_{i+j+1}w_{i+j}.
\end{equation*}
The claim then follows by induction.  To complete the proof, we note that $\rho(F^{\infty})<1$ \cite{bacsar2008h} and invoke Lemma \ref{the:lyapunov} to replace $P_{i}^{\infty}$ and $P_{i+1}^{\infty}$ by $P^{\infty}$ in \eqref{eq:cost_to_go_infty}. We note that while $P^{\infty}$, $v_i^{\infty}$ and $q_i^{\infty}$ are independent of $x_i$, the last two depend on the noise realisation $w$.  $\hfill \square$

\textit{Remark}: Lemma \ref{the:induction} also holds for any stabilising state feedback matrix $K^s$, with the coefficients in the cost-to-go appropriately defined.

The optimal offline controller, as defined in section \ref{sec:regret_definition}, has access to all future disturbances $w$ and minimizes the cost function \eqref{eq:cost_T} without constraining the inputs to a policy set. It is shown in \cite{goel2020power} that for the infinite horizon case the controller has the following form
\begin{equation}
    u_t^{\star} = -Kx_t -(R+B^{\top}PB)^{-1}B^{\top}\sum_{i=0}^{\infty}(F^{\top})^iPw_{t+i},
    \label{eq:optimal_offline_infinite}
\end{equation}
where $K = -(R+B^{\top}PB)^{-1}B^{\top}PA$, $P$ is the solution of the discrete ARE and $F:= A - BK$, with $\rho(F)<1$. Moreover, the cost-to-go  at a state $x_i$, $J(u^{\star},w;x_i)$ has the same extended quadratic structure as in \eqref{eq:cost_to_go_infty}  with the following coefficients
\begin{align}
P & = F^{ \top}PF + Q + K^{\top}RK,
\label{eq:P_t}\\
v_i & = F^{  \top}(2Pw_i + v_{i+1}) = 2 \sum_{j=0}^{\infty}(F^{  \top})^{j+1}Pw_{i+j},\label{eq:v_i_optimal_offline}\\
q_i & = q_{i+1} + w_i^{\top}v_{i+1} + w_i^{\top}Pw_i - {G_i}^{\top}H{G_i},
\end{align}
where ${G_i}: = \sum_{j=0}^{\infty}(F^{ \top})^{j}Pw_{i+j}$ and $H:= B(R+B^{\top}PB)^{-1}B^{\top}$ \cite{goel2020power, yu2020power}. Using the result of Lemma \ref{the:induction}, the  regret of the $\mathcal{H}_{\infty}$  controller for a given $w$ is then
\begin{equation}
    \begin{split}
    & \mathcal{R}(\mathcal H_{\infty},w) := \lim_{T \rightarrow \infty}\mathcal{R}_T(\mathcal H_{\infty},w)\\
    &=J(u^{\infty},w;x_0) - J(u^{\star},w;x_0) \\
    &= x_0^{\top}(P^{\infty} - P)x_0 + x_0^{\top}(v_0^{\infty}-v_0)+q_0^{\infty} - q_0.
\end{split}
\label{eq:regret_Hinf}
\end{equation}

The \textit{disturbance-reality gap} $\Delta w$ is defined as  
\begin{equation}
   \Delta w : =  w - w^{\star}.
   \label{eq:disturbance_reality_gap}
\end{equation}
This vector is the difference between the disturbance $w$, experienced by the system  and the worst-case disturbance,  $w^{\star}$, assumed by the $\mathcal{H}_{\infty}$ controller,  defined in \eqref{eq:w_worst_inf}. The main result of this paper is formulated in the following theorem.
\begin{theorem}
The $\mathcal{H}_{\infty}$ controller, that solves the disturbance attenuation problem \eqref{eq:h_inf_problem} attains dynamic regret
\begin{equation*}
    \mathcal{R}(\mathcal{H}_{\infty}, w) \leq  k_1\|\Delta w\| + k_2\|\Delta w\|^2
\end{equation*}
for all initial states in \eqref{eq:x0_set} and constants $k_1,k_2 \in \mathbb{R}^+$, given below in \eqref{eq:k_1_and_k_2}.
\end{theorem}

\textit{Proof}:  Regret \eqref{eq:regret_Hinf} can be written in terms of $\Delta w$ and $w^{\star}$
\begin{align*}
\begin{split}
     &\mathcal{R}(\mathcal{H}_{\infty},w) = x_0^{\top}(P^{\infty} - P)x_0\\
     &+\sum_{i=0}^{\infty}\biggl(2x_0^{\top} \left((F^{\infty  \top})^{i+1}P^{\infty} - (F^{\top})^{i+1}P \right)(w^{ \star}_i + \Delta w_i) \\
     &+ (w^{ \star}_i+ \Delta w_i)^{\top}(v_{i+1}^{\infty} - v_{i+1}) +{G_i}^{\top}H{G_i}\\
     &+ (w^{ \star}_i + \Delta w_i)^{\top}(P^{\infty} - P)(w^{ \star}_i + \Delta w_i)\biggr).
     \end{split}
\end{align*}
Since the $\mathcal{H}_{\infty}$ controller is in saddle point equilibrium with $w^{\star}$,  the optimal offline controller with  knowledge of the future disturbances, will attain the same cost as the $\mathcal{H}_{\infty}$ controller if $\Delta w = 0$. Hence $R(\mathcal{H}_\infty,w^*)=0$, leaving
\begin{align*}
\begin{split}
     &\mathcal{R}(\mathcal{H}_\infty,w) =\\
     &\lim_{T \rightarrow \infty }\sum_{i=0}^{T-1}\biggl(\underbrace{ 2x_0^{\top}\left((F^{\infty  \top})^{i+1}P^{\infty} - (F^{  \top})^{i+1}P \right)\Delta w_i}_\text{$a_i$}  \\
     &+ \underbrace{w^{\star \top}_i(v_{i+1}^{\infty \Delta w} - v_{i+1}^{\Delta w})}_\text{$b_i$}  +\underbrace{\Delta w^{\top}_i(v_{i+1}^{\infty w} - v_{i+1}^w)}_\text{$c_i$}\\
     &+ \underbrace{\Delta w_i^{\top}(P^{\infty} - P)\Delta w_i}_\text{$d_i$} + \underbrace{2\Delta w_i^{\top}(P^{\infty} - P)w^{ \star}_i}_\text{$e_i$}\\
     &+ \underbrace{{G_i}^{\Delta w \top}H{G_i}^{\Delta w}}_\text{$f_i$}+ \underbrace{2{G_i}^{\Delta w \top}H{G_i}^{w^{\star}}}_\text{$g_i$}\biggr),
     \end{split}
\end{align*}
where $v_i^{ \Delta w} := F^{  \top}(2P \Delta w + v^{\Delta w}_{i+1})$, $v_i^{\infty \Delta w}: = F^{ \infty \top}(2P^{\infty} \Delta w + v^{\infty \Delta w}_{i+1})$ and ${G_i}^{\Delta w}: =\sum_{j=0}^{T-i-1}(F^{ \top})^{j}P \Delta w_{i+j}$; the corresponding expressions with $w^{\star}$ and $w$ are defined analogously. This reformulation of regret is then used to upper bound it in terms of the norm of $\Delta w$. For $d_i$ and $e_i$
\begin{align*}
    \lim_{T\rightarrow \infty}\sum_{i=0}^{T-1}d_i &= \Delta w ^{\top}\biggl(\Delta P \otimes I_{nT} \biggr) \Delta w \leq 2 \|\Delta w \|^2 \|\bar{P}\| \\
    \lim_{T\rightarrow \infty}\sum_{i=0}^{T-1} e_i &= 2 \Delta w ^{\top}\biggl(\Delta P\otimes I_{nT} \biggr)   w^{\star} \leq 4 \|\Delta w \| \|\bar{P}\|,
\end{align*}
where $\Delta P:=P^{\infty}-P$, $\|\bar{P}\| := \max\{\|P^{\infty}\|,\|P\|\}$, and using  $\|w^{ \star} \| = 1$ and the fact \cite{lancaster1972norms} that $\|A \otimes B\| = \|A\| \|B\|$. The sum of terms $a_i$ can be written as
\begin{equation*}
    \sum_{i=0}^{T-1}a_i = 2x_0^{\top} \Delta L_v \Delta w,
\end{equation*}
where $\Delta L_v \in \mathbb{R}^{n \times nT}$ is a block matrix with the term $(F^{\infty  \top})^{i}P^{\infty} - (F^{  \top})^{i}P$ on its $i$-th block column for all $1\leq i\leq T$. From Gelfand's formula it can be shown that there exists a constant $c>1$ such that $\|F^i\|\leq c\lambda^i$ and $\|(F^{\infty})^i\|\leq c(\lambda^{\infty})^i$ for all $i\geq1$ with $\lambda: = \frac{1+\rho(F)}{2}<1$, $\lambda^{\infty} = \frac{1+\rho(F^{\infty})}{2}$ , since $\rho(F) < 1$, $\rho(F^{\infty}) < 1$ . Hence
\begin{align*}
\begin{split}
  \lim_{T\rightarrow \infty} \sum_{i=0}^{T-1} a_i &\leq 2c\|x_0\|\|\bar{P}\|\left(\frac{\lambda^{\infty}}{1-\lambda^{\infty}} + \frac{\lambda}{1-\lambda}\right)\|\Delta w\| \\
    & \leq 4cX\|\Delta w\|\|\bar{P}\|\biggl(\frac{\bar{\lambda}}{1-\bar{\lambda}}\biggr),
    \end{split}
\end{align*}
$\bar{\lambda}:= \max\{\lambda^{\infty},\lambda\}$. For the term with $c_i$
\begin{equation*}
    \sum_{i=0}^{T-1}c_i = 2 \Delta w_{[0:T-2]}^{\top}\Delta L_u w_{[1:T-1]},
\end{equation*}
where $\Delta L_u \in \mathbb R^{n(T-1) \times n(T-1)}$ is an upper triangular block Toeplitz matrix, such that for all $1 \leq i <T$ and $i \leq j < T$, the matrix on the $i$-th block row and $j$-th block column is $(F^{\infty  \top})^{j-i+1}P^{\infty} - (F^{\top})^{j-i+1}P$. It follows that
\begin{align*}
\begin{split}
     \lim_{T\rightarrow \infty}\sum_{i=0}^{T-1}c_i &\leq 2c\|\Delta w\|\left(\|P^{\infty}\|\frac{\lambda^{\infty}}{1-\lambda^{\infty}} + \|P\|\frac{\lambda}{1-\lambda}\right)\\
      &\leq 4c\|\Delta w\|\|\bar{P}\|\biggl(\frac{\bar{\lambda}}{1-\bar{\lambda}}\biggr),
     \end{split}
\end{align*}
where we have used the properties of block Toeplitz matrices. We can similarly get the same bound for $\sum_{i=0}^{T-1}b_i$. The sum of terms $f_i$ and $g_i$ can be written as
\begin{align*}
     \sum_{i=0}^{T-1}f_i&=\Delta w^{\top}L_g^{\top}(H \otimes I_{nT})L_g\Delta w,\\
     \sum_{i=0}^{T-1}g_i &=2\Delta w^{\top}L_g^{\top}(H \otimes I_{nT})L_g w,
\end{align*}
where  $L_g \in \mathbb{R}^{nT}$ is an upper triangular block Toeplitz matrix, such that for all $1 \leq i \leq T$ and $i \leq j \leq T$, the matrix on the $i$-th block row and $j$-th block column is $(F^{\top})^{j-i}P$. Upper bounds for both can then be similarly found
\begin{align*}
    \lim_{T\rightarrow \infty} \sum_{i=0}^{T-1}f_i \leq \|\Delta w\|^2 \|H\|\|P\|^2\frac{c^2}{(1-\lambda)^2},\\
    \lim_{T\rightarrow \infty} \sum_{i=0}^{T-1}g_i \leq2 \|\Delta w\|\|H\|\|P\|^2\frac{c^2}{(1-\lambda)^2}.
\end{align*}
Summing the terms and setting
\begin{align}
\begin{split}
    k_2 &= 2\|\bar{P}\|+ \|H\|\|P\|^2\frac{c^2}{(1-\lambda)^2}\\
    k_1 &= 4\|\bar{P}\|+ 4 c\|\bar{P}\|\left(2+X\right)\biggl(\frac{\bar{\lambda}}{1-\bar{\lambda}}\biggr)\\
    &+2 \|H\|\|P\|^2\frac{c^2}{(1-\lambda)^2},
    \end{split}
    \label{eq:k_1_and_k_2}
\end{align}
completes the proof. $\hfill \square$

We note that the constraint of the noise signal having a unit energy is without loss of generality and the same result can also be attained by modifying the set \eqref{eq:x0_set}. \\

\subsubsection{Finite Horizon}
A similar bound is obtained for the dynamic regret of the finite horizon controller.
\begin{theorem}
The $\mathcal{H}_{\infty}$ controller, that solves the disturbance attenuation problem \eqref{eq:h_inf_problem} for a horizon length $T$,  attains dynamic regret
\begin{equation*}
    \mathcal{R}_T(\mathcal{H}_{\infty}, w) \leq k'_1\|\Delta w\| + k'_2\|\Delta w\|^2
\end{equation*}
for all initial states in \eqref{eq:x0_set} and constants $k'_1,k'_2 \in \mathbb{R}^+$, given below in \eqref{eq:k_1_and_k_2_finite}.
\end{theorem}
\textit{Proof}: The proof follows closely the structure for the infinite horizon case. Using the finite horizon $\mathcal{H}_{\infty}$ controller \eqref{eq:control_strategy} and following the induction arguments in Lemma \ref{the:induction} a similar extended quadratic expression for the cost-to-go is achieved. Specifically, for this controller the cost-to-go at time step $i$, $0 \leq i <T$ is given as $J_i^T(u^{\infty},w,x_i)=x_i^{\top}P^{\infty}_ix_i + x_i^{\top}v^{\infty}_i +q^{\infty}_i$, where
\begin{align*}
P^{\infty}_i & = F^{\infty \top}_i P^{\infty}_{i+1}F_i^{\infty} + Q + K^{\infty \top}_iRK_i^{\infty}
\label {eq:Pg_t}\\
v^{\infty}_i & = 2 \sum_{j=0}^{T-i-1}\Phi^{\infty}(i+j+1,i)^{\top}P_{i+j+1}^{\infty}w_{i+j}\\
q^{\infty}_i & = q^{\infty}_{i+1} + w^{\top}_iv^{\infty}_{i+1} + w^{\top}_iP^{\infty}_{i+1}w_i,
\end{align*}
where the state transition matrix is defined as
\begin{equation}
    \Phi^{\infty}(t,t_0): = \begin{cases}
    F^{\infty}_{t-1} F^{\infty}_{t-2}\dots F^{\infty}_{t_0}, & \text{if $t>t_0$}.\\
    I_n, & \text{$t=t_0$}.
    \label{eq:phi_inf}
  \end{cases}
\end{equation}

The optimal offline controller for the finite horizon case can be written in the following form \cite{goel2020power, zhang2021regret}
\begin{equation*}
    u_t^{\star} = -K_tx_t-\sum_{i=0}^{T-t-1}K_{t,i}^w w_{t+i},
\end{equation*}
where for all $0\leq i<T-t-1$,
\begin{equation*}
    K_{t,i}^w = (R+B^{\top}P_{t+1}B)^{-1}B^{\top}\Phi(t+i+1,t+1)^{\top}P_{t+i+1}.
\end{equation*}
Here $P_t$ is the solution of the difference Ricatti equation arising from the standard LQR problem, and $K_t= (R+B^{\top}P_{i+1}B)^{-1}B^{\top}P_{t+1}A$ is the associated optimal gain, both at time $0\leq t < T$. The state transition matrix $\Phi$ is defined analogously to \eqref{eq:phi_inf} with the corresponding $F_t$-s. The cost-to-go is again represented as an extended quadratic function $J(u^{\star}, w,x_i) = x_i^{\top}P_ix_i + x_i^{\top}v_i +q_i$, where
\begin{align*}
P_i & = F^{\top}_i P_{i+1}F_i + Q + K^{\top}_iRK_i\\
v_i & = 2 \sum_{j=0}^{T-i-1}\Phi(i+j+1,i)^{\top}P_{i+j+1}w_{i+j}\\
q_i & = q_{i+1} + w^{\top}_iv_{i+1} + w^{\top}_iP_{i+1}w_i- \underline{G_i}^{\top}H_i\underline{G_i},
\end{align*}
where $\underline{G_i}: = \sum_{j=0}^{T-i-1}\Phi(i+j+1,i+1)^{\top}P_{i+j+1}w_{i+j}$ and $H_i:= B(R+B^{\top}P_{i+1}B)^{-1}B^{\top}$.

The regret of the finite horizon $\mathcal{H}_{\infty}$ controller is then the difference of the two extended quadratic functions,
\begin{equation*}
\mathcal{R}_T(\mathcal H_{\infty}) =x_0^{\top}(P^{\infty}_0 - P_0)x_0 + x_0^{\top}(v_0^{\infty}-v_0)+q_0^{\infty} - q_0.
\end{equation*}
Substituting the expressions for the coefficients, the above can be written in terms of the \textit{disturbance-reality gap}. Using the same argument of equal costs for the worst-case disturbance signal,  the following terms are left 
\begin{align*}
\begin{split}
     &\mathcal{R}_T(\mathcal{H}_{\infty},w) =\sum_{i=0}^{T-1}\biggl(\underbrace{\underline{G_i}^{\Delta w \top}H_i\underline{G_i}^{\Delta w}}_\text{$f_i$}+\underbrace{2\underline{G_i}^{\Delta w \top}H_i\underline{G_i}^{w^{\star}}}_\text{$g_i$} \\
     &+\underbrace{2x_0^{\top} \left(\Phi^{\infty}(i+1,0)^{\top}P_{i+1}^{\infty} - \Phi(i+1,0)^{\top}P_{i+1} \right)\Delta w_i}_\text{$a_i$}  \\
     &+\underbrace{w^{\star \top}_i(v_{i+1}^{\infty \Delta w} - v_{i+1}^{\Delta w})}_\text{$b_i$}  + \underbrace{\Delta w^{\top}_i(v_{i+1}^{\infty w} - v_{i+1}^w)}_\text{$c_i$}\\
     &+ \underbrace{\Delta w_i^{\top}(P^{\infty}_{i+1} - P_{i+1})\Delta w_i}_\text{$d_i$} +\underbrace{2\Delta w_i^{\top}(P^{\infty}_{i+1} - P_{i+1})w^{ \star}_i}_\text{$e_i$} \biggr),\\
     \end{split}
\end{align*}
where $v_i^{ \Delta w} := F_i^{  \top}(2P_{i+1} \Delta w + v^{\Delta w}_{i+1})$, $v_i^{\infty \Delta w}: = F_i^{ \infty \top}(2P_{i+1}^{\infty} \Delta w + v^{\infty \Delta w}_{i+1})$ and ${G_i}^{\Delta w}: =\sum_{j=0}^{T-i-1}\Phi(i+j+1,i+1)^{\top}P_{i+j+1} \Delta w_{i+j}$; the corresponding expressions with $w^{\star}$ and $w$ are defined analogously. Defining  $\bar{P}'\succeq P_i, P_i^{\infty}$ $\quad \forall 0 \leq i \leq T$ the following bounds can then be achieved
\begin{align*}
  \sum_{i=0}^{T-1}d_i & \leq 2 \|\Delta w \|^2 \|\bar{P}'\| \\
\sum_{i=0}^{T-1} e_i & \leq 4 \|\Delta w \| \|\bar{P}'\|.
\end{align*}
The sum of terms $a_i$ can be written as
\begin{equation*}
    \sum_{i=0}^{T-1}a_i = 2x_0^{\top} \Delta L^f_v \Delta w,
\end{equation*}
where $\Delta L^f_v \in \mathbb{R}^{n \times nT}$ is a block matrix with the term $\Phi^{\infty}(i,0)^{\top}P_{i}^{\infty} - \Phi(i,0)^{\top}P_{i}$ on its $i$-th block column for all $1\leq i\leq T$. Using the results of exponential stability for finite horizon LQR \cite{zhang2021regret} and defining $\bar{\tau} := \sqrt{\frac{\lambda_{max}(\bar{P}')}{\lambda_{min}(Q)}}$ and $\bar{\eta} := \sqrt{1 - \frac{1}{\bar{\tau}^2}}<1$, the following bound is achieved
\begin{equation*}
 \sum_{i=0}^{T-1} a_i \leq 4X \|\bar{P}'\| \bar{\tau}\bar{\eta}\left(\frac{1-\bar{\eta}^{T}}{1-\bar{\eta}}\right)\|\Delta w\|.
 \end{equation*}
For the sum of terms $c_i$
\begin{equation*}
    \sum_{i=0}^{T-1}c_i = 2 \Delta w_{[0:T-2]}^{\top}\Delta L_u^f w_{[1:T-1]},
\end{equation*}
where $\Delta L_u^f \in \mathbb R^{n(T-1) \times n(T-1)}$ is an upper triangular block Toeplitz matrix, such that for all $1 \leq i <T$ and $i \leq j < T$, the matrix on the $i$-th block row and $j$-th block column is $\Phi^{\infty}(j+1,i)^{\top}P_{j+1}^{\infty} - \Phi(j+1,i)^{\top}P_{j+1}$. It follows that
\begin{equation*}
 \sum_{i=0}^{T-1}c_i \leq 4 \|\Delta w\|\|\bar{P}'\|\bar{\tau}\bar{\eta}\left(\frac{1-\bar{\eta}^{T-1}}{1-\bar{\eta}}\right)
\end{equation*}
We can similarly get the same bound for $\sum_{i=0}^{T-1}b_i$. The sum of last two terms $f_i$ and $g_i$ can be written as
\begin{align*}
     \sum_{i=0}^{T-1}f_i&=\Delta w^{\top}L^{f \top}_g H^d L^{f}_g\Delta w,\\
     \sum_{i=0}^{T-1}g_i &=2\Delta w^{\top}L^{f \top}_g H^d L^{f}_g w,
\end{align*}
where  $L_g^f \in \mathbb{R}^{nT\times nT}$ is an upper triangular block Toeplitz matrix, such that for all $1 \leq i \leq T$ and $i \leq j \leq T$, the matrix on the $i$-th block row and $j$-th block column is $\Phi(j,i)^{\top}P_{j}$ and $H^d\in \mathbb{R}^{nT \times nT}$ is a block diagonal matrix with $H_{i}$ on its $i$-th block diagonal entry. Defining  $\bar{H}\succeq B(R+B^{\top}P_iB)^{-1}B^{\top} \quad  \forall \quad 0 < i \leq T$  upper bounds for both terms are then given as
\begin{align*}
    \sum_{i=0}^{T-1}f_i \leq\bar{\tau}^2 \|\Delta w\|^2 \|\bar{H}\|\|\bar{P}'\|^2\frac{(1-\bar{\eta}^{T})^2}{(1-\bar{\eta})^2},\\
    \sum_{i=0}^{T-1}g_i \leq2\bar{\tau}^2 \|\Delta w\|\|\bar{H}\|\|\bar{P}'\|^2\frac{(1-\bar{\eta}^{T})^2}{(1-\bar{\eta})^2}.
\end{align*}
Summing the terms and setting the constants $k_1'$ and $k_2'$ as follows
\begin{align}
\begin{split}
    k_2' &= \|\bar{P}'\| \biggl(2+\bar{\tau}^2\|\bar{H}\|\|\bar{P}'\|\frac{(1-\bar{\eta}^{T})^2}{(1-\bar{\eta})^2}\biggr)\\
    k_1' &= 2\|\bar{P}'\|\biggl(2+ 2\bar{\tau}\bar{\eta}(2+X)\left(\frac{1-\bar{\eta}^{T}}{1-\bar{\eta}}\right)\\
        &+\bar{\tau}^2\|\bar{H}\|\|\bar{P}'\|\frac{(1-\bar{\eta}^{T})^2}{(1-\bar{\eta})^2}\biggr),
    \end{split}
    \label{eq:k_1_and_k_2_finite}
\end{align}
completes the proof. $\hfill \square$
\section{CE OPTIMISTIC CONTROLLER}
In this section, the certainty equivalent  optimistic controller that has an inaccurate prediction $\bar{w} \in \mathcal{L}_2$ of the disturbance signal and acts optimally with respect to it is considered. It solves the optimization problem \eqref{eq:cost_T} subject to $x_{t+1}= Ax_t+Bu_t+\bar{w}_t$. For simplicity only  the infinite horizon case is considered, however, the results for the finite horizon are derived analogously. The infinite horizon CE controller is the same as in \eqref{eq:optimal_offline_infinite}, only with feedback on $\bar{w}$
\begin{equation}
    u_t^{\mbox{CE}} = -Kx_t -(R+B^{\top}PB)^{-1}B^{\top}\sum_{i=0}^{\infty}F^{\top i}P\bar{w}_{t+i}.
    \label{eq:CE_controller}
\end{equation}
The dynamic regret for this controller is shown to be proportional to $\| \Delta \bar{w}\|^2$, where $\Delta \bar{w}:=  w -\bar{w} $, is the error vector between the predicted and the observed true disturbance. The result is formulated in the following proposition.
\begin{proposition}
The CE optimistic controller \eqref{eq:CE_controller} attains dynamic regret that is upper bounded by
\begin{equation*}
    \mathcal{R}(\mbox{CE},w) \leq\|\Delta \bar{w}\|^2 \|H\|\|P\|^2\frac{c^2}{(1-\lambda)^2}.
\end{equation*}
\end{proposition}
\textit{Proof:} We start with the same induction hypothesis that the cost-to-go at timestep $i\geq0$ is $J(u^{CE},w;x_i) = x_i^{\top}P_ix_i + \bar{v}_i^{\top}x_i +\bar{q}_i$. For the timestep $T$ we have trivially $P_T =Q_T, \bar{v}_T =0$ and $\bar{v}_T = 0$. We note that the state feedback matrix $K$ of $u_t^{\mbox{CE}}$ is stabilising in this case as well, and follow the same technique as in Lemma \ref{the:induction} in the limit of $T \rightarrow \infty$ to get 
\begin{align*}
\begin{split}
   &J(u^{\mbox{CE}},w;x_i) = \bar{q}_{i+1} + x_i^{\top}(Q+A^{\top}PA - A^{\top}PHPA)x_i \\
   &+x_i^{\top}(F^{\top}\bar{v}_{i+1}+2F^{\top}Pw_i)+ w_i^{\top}Pw_i+w_i^{\top}\bar{v}_{i+1}\\
   &+(\frac{\bar{v}_{i+1}}{2} - \sum_{j=1}^{\infty}(F^{\top})^jP\bar{w}_{i+j})^{\top}H(\frac{\bar{v}_{i+1}}{2}-\sum_{j=1}^{\infty}(F^{\top})^jP\bar{w}_{i+j})\\
   &-(P\bar{w}_i + \frac{\bar{v}_{i+1}}{2})^{\top}H(P\bar{w}_i + \frac{\bar{v}_{i+1}}{2})\\
   &-2\biggl(\sum_{j=0}^{\infty}(F^{\top})^jP\bar{w}_{i+j}\biggr)^{\top}HP\Delta \bar{w}_i.
    \end{split}
\end{align*}
From the above, it can be concluded, that in order for the induction hypothesis to hold, $P$ needs to be the solution of the associated ARE for the problem, $\bar{v}_i = v_i $ is the same as for the optimal offline controller  \eqref{eq:v_i_optimal_offline}, and,
\begin{align*}
\begin{split}
   &\bar{q}_i = \bar{q}_{i+1} +(\sum_{j=1}^{\infty}(F^{\top})^jP\Delta \bar{w}_{i+j})^{\top}H(\sum_{j=1}^{\infty}(F^{\top})^jP\Delta \bar{w}_{i+j})\\
   &-2\biggl(\sum_{j=0}^{\infty}(F^{\top})^jP\bar{w}_{i+j}\biggr)^{\top}HP\Delta w_i+ w_i^{\top}Pw_i+w_i^{\top}\bar{v}_{i+1}\\
   &-\bar{G}_i^{\top}H\bar{G}_i,\\
    \end{split}
\end{align*}
where $\bar{G}_i := -P\Delta \bar{w}_i + \sum_{j=0}^{\infty}(F^{\top})^jP{w}_{i+j}$.
The regret of the algorithm then equals to the difference between the constant terms,
\begin{align*}
\begin{split}
&\mathcal{R}(\mbox{CE},w) =J(u^{\mbox{CE}},w;x_0) - J(u^{\star},w,x_0)= \bar{q}_0 - q_0\\
& = \sum_{i=0}^{\infty}\biggl(\sum_{j=0}^{\infty}(F^{\top})^jP\Delta \bar{w}_{i+j}\biggr)^{\top}H\biggl(\sum_{j=0}^{\infty}(F^{\top})^jP\Delta \bar{w}_{i+j}\biggr).
\end{split}
\end{align*}
Similar to the proof of the $\mathcal{H}_\infty$ controller, the required upper bound is then achieved.
\hfill $\square$

Comparing the above with the regret for the infinite horizon $\mathcal{H}_{\infty}$ controller we note that in addition to the dependence on $\|\Delta w\|$, the regret bound for the $\mathcal{H}_{\infty}$ also has additional terms in the coefficient of $\|\Delta w\|^2$. Thus, given  equal  $\|\Delta w\|$ and $\|\Delta \bar{w}\|$, $\mathcal{H}_{\infty}$ has a strictly higher regret upper bound compared to the certainty equivalent controller. The additional terms in the upper bound of $\mathcal{H}_{\infty}$ regret are due to the ``sub-optimal'' state feedback gain on the state. This results in a mismatch between $P^{\infty}$ and $P$, as well as $v_0^{\infty}$ and $v_0$, the coefficients of the initial state; this is not the case for the CE controller. This makes the \textit{gap} explicitly dependent on the initial state leading to additional terms in the regret bound.
\section{NUMERICAL EXAMPLE}
\begin{figure}
     \center
         \includegraphics[width=\linewidth]{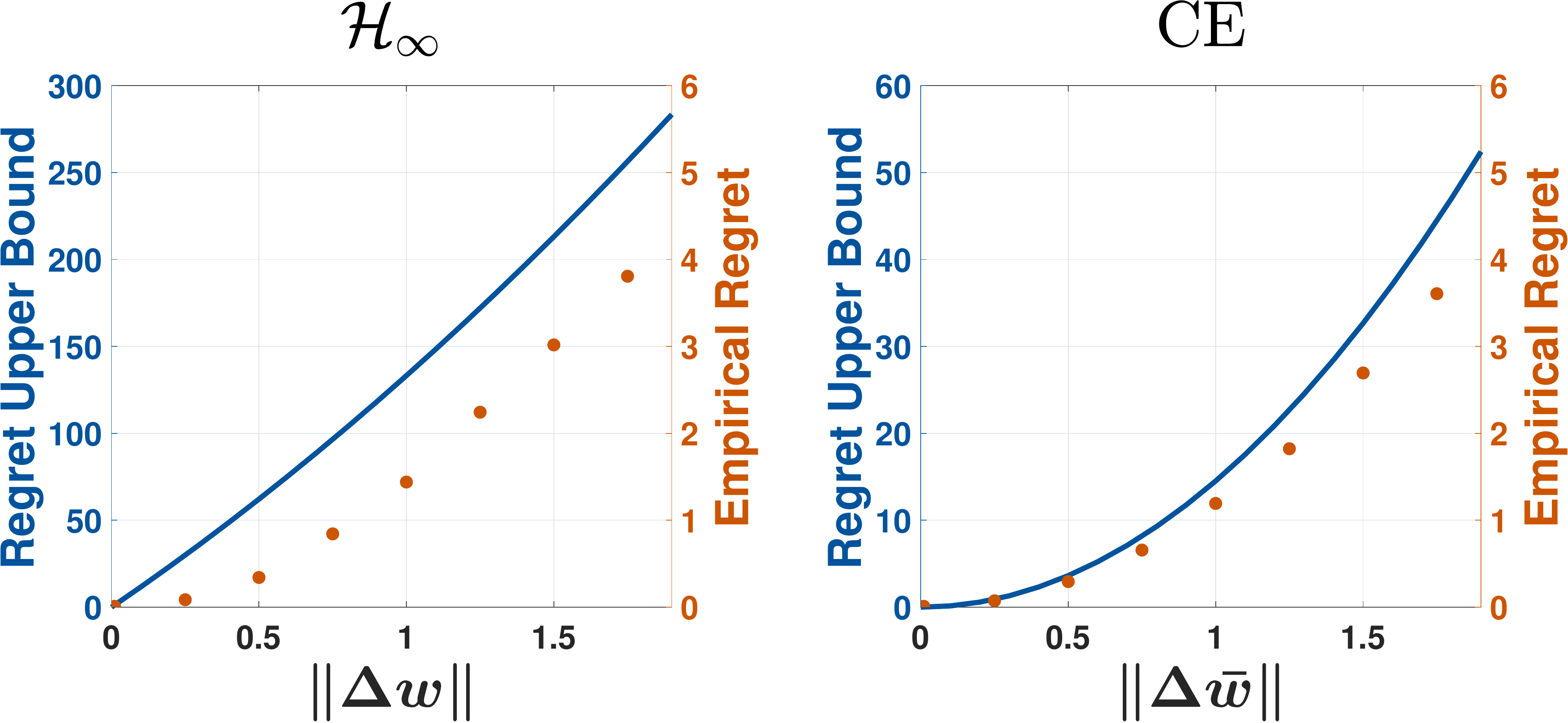}
         \caption{The theoretical regret upper bounds and worst-case simulated regret for the $\mathcal{H}_\infty$ and $\mbox{CE}$ controllers.  }
         \label{fig:numerics}
\end{figure}

A simple system, also considered in \cite{didinsky1992design}, with $A=1, B=1$ and cost matrices $Q=1, R=1$ is controlled using both the $\mathcal{H}_{\infty}$ and the $\mbox{CE}$ finite horizon controllers with $T=100$.  For each fixed $\|\Delta w\|$ and $\|\Delta \bar{w}\|$, a number of random noise signals are generated. The system evolution is then simulated starting from an initial condition $x_0 = 4 \in X_s$. The parameter $\bar{\gamma}$ is found for this initial state using a trial-and-error method as described in \cite{didinsky1992design}. The dynamic regret of both controllers is then calculated and the highest regret (for each $\|\Delta w\|$ and $\|\Delta \bar{w}\|$) is plotted in Figure \ref{fig:numerics} along with the upper bounds obtained in this work. It is inferred from the plots that the analytic bounds capture the order of the empirically calculated worst-case regret. Preliminary numerical tests show that they can become tighter for certain adversarial noise signals.

\section{CONCLUSIONS}
The $\mathcal{H}_{\infty}$ algorithm is considered in the context of dynamic regret to characterize its extra cost due to planning for the worst-case disturbance realization. The upper bound of this regret is shown to scale with the norm of the \textit{gap} between the worst-case predicted by the $\mathcal{H}_\infty$ controller and the true one. This result is compared with the CE optimistic controller with erroneous predictions. While both controllers have similar regret upper bound structures, for equal \textit{disturbance-reality gap} and \textit{prediction error} norms, $\mathcal{H}_{\infty}$'s regret attains a strictly higher upper bound. A  numerical example is presented to show that the order of the worst-case simulated regret is captured by the theoretical bounds. A possible direction for further research, is the consideration of the case where no saddle point exists for the $\mathcal{H}_{\infty}$ problem. The results can also provide insights on the development of algorithms that estimate the future disturbances online.

\textit{Acknowledgements}: The authors thank Efe Balta for fruitful discussions on the topic.
\bibliographystyle{ieeetr}
\bibliography{main.bib}

\begin{thebibliography}{10}

\bibitem{hassibi1999indefinite}
B.~Hassibi, A.~H. Sayed, and T.~Kailath, {\em Indefinite-Quadratic estimation
  and control: a unified approach to H 2 and H$\infty$ theories}.
\newblock SIAM, 1999.

\bibitem{zhou1998essentials}
K.~Zhou and J.~C. Doyle, {\em Essentials of robust control}, vol.~104.
\newblock Prentice hall Upper Saddle River, NJ, 1998.

\bibitem{green2012linear}
M.~Green and D.~J. Limebeer, {\em Linear robust control}.
\newblock Courier Corporation, 2012.

\bibitem{books/daglib/0016248}
N.~Cesa-Bianchi and G.~Lugosi, {\em Prediction, learning, and games.}
\newblock Cambridge University Press, 2006.

\bibitem{hazan2020nonstochastic}
E.~Hazan, S.~Kakade, and K.~Singh, ``The nonstochastic control problem,'' in
  {\em Algorithmic Learning Theory}, pp.~408--421, PMLR, 2020.

\bibitem{goel2020power}
G.~Goel and B.~Hassibi, ``The power of linear controllers in {LQR} control,''
  {\em arXiv preprint arXiv:2002.02574}, 2020.

\bibitem{dean2018regret}
S.~Dean, H.~Mania, N.~Matni, B.~Recht, and S.~Tu, ``Regret bounds for robust
  adaptive control of the linear quadratic regulator,'' {\em Advances in Neural
  Information Processing Systems}, vol.~31, 2018.

\bibitem{simchowitz2020naive}
M.~Simchowitz and D.~Foster, ``Naive exploration is optimal for online {LQR},''
  in {\em International Conference on Machine Learning}, pp.~8937--8948, PMLR,
  2020.

\bibitem{agarwal2019online}
N.~Agarwal, B.~Bullins, E.~Hazan, S.~Kakade, and K.~Singh, ``Online control
  with adversarial disturbances,'' in {\em International Conference on Machine
  Learning}, pp.~111--119, PMLR, 2019.

\bibitem{foster2020logarithmic}
D.~Foster and M.~Simchowitz, ``Logarithmic regret for adversarial online
  control,'' in {\em International Conference on Machine Learning},
  pp.~3211--3221, PMLR, 2020.

\bibitem{yu2020power}
C.~Yu, G.~Shi, S.-J. Chung, Y.~Yue, and A.~Wierman, ``The power of predictions
  in online control,'' {\em arXiv preprint arXiv:2006.07569}, 2020.

\bibitem{zhang2021regret}
R.~Zhang, Y.~Li, and N.~Li, ``On the regret analysis of online {LQR} control
  with predictions,'' in {\em 2021 American Control Conference (ACC)},
  pp.~697--703, IEEE, 2021.

\bibitem{muthirayan2021robust}
D.~Muthirayan, D.~Kalathil, and P.~P. Khargonekar, ``Online robust control of
  linear dynamical systems with prediction,'' {\em arXiv preprint
  arXiv:2111.15063}, 2021.

\bibitem{goel2020regret}
G.~Goel and B.~Hassibi, ``Regret-optimal control in dynamic environments,''
  {\em arXiv preprint arXiv:2010.10473}, 2020.

\bibitem{sabag2021regret}
O.~Sabag, G.~Goel, S.~Lale, and B.~Hassibi, ``Regret-optimal controller for the
  full-information problem,'' in {\em 2021 American Control Conference (ACC)},
  pp.~4777--4782, IEEE, 2021.

\bibitem{martin2022safe}
A.~Martin, L.~Furieri, F.~D{\"o}rfler, J.~Lygeros, and G.~F. Trecate, ``Safe
  control with minimal regret,'' {\em arXiv preprint arXiv:2203.00358}, 2022.

\bibitem{didier2022system}
A.~Didier, J.~Sieber, and M.~N. Zeilinger, ``A system level approach to regret
  optimal control,'' {\em arXiv preprint arXiv:2202.13763}, 2022.

\bibitem{bacsar2008h}
T.~Ba{\c{s}}ar and P.~Bernhard, {\em H-infinity optimal control and related
  minimax design problems: a dynamic game approach}.
\newblock Springer Science \& Business Media, 2008.

\bibitem{muthirayan2021online}
D.~Muthirayan, J.~Yuan, D.~Kalathil, and P.~P. Khargonekar, ``Online learning
  for receding horizon control with provable regret guarantees,'' {\em arXiv
  preprint arXiv:2111.15041}, 2021.

\bibitem{hazan2019introduction}
E.~Hazan, ``Introduction to online convex optimization,'' {\em arXiv preprint
  arXiv:1909.05207}, 2019.

\bibitem{bacsar1998dynamic}
T.~Ba{\c{s}}ar and G.~J. Olsder, {\em Dynamic noncooperative game theory}.
\newblock SIAM, 1998.

\bibitem{limebeer1989discrete}
D.~Limebeer, M.~Green, and D.~Walker, ``Discrete-time $\mathcal{H}_\infty$
  control,'' in {\em Proceedings of the 28th IEEE Conference on Decision and
  Control,}, pp.~392--396, IEEE, 1989.

\bibitem{didinsky1992design}
G.~Didinsky and T.~Ba{\c{s}}ar, ``Design of minimax controllers for linear
  systems with non-zero initial states under specified information
  structures,'' {\em International Journal of Robust and Nonlinear Control},
  vol.~2, no.~1, pp.~1--30, 1992.

\bibitem{mageirou1976values}
E.~Mageirou, ``Values and strategies for infinite time linear quadratic
  games,'' {\em IEEE Transactions on Automatic Control}, vol.~21, no.~4,
  pp.~547--550, 1976.

\bibitem{willems1971least}
J.~Willems, ``Least squares stationary optimal control and the algebraic
  riccati equation,'' {\em IEEE Transactions on automatic control}, vol.~16,
  no.~6, pp.~621--634, 1971.

\bibitem{dunford1963spectral}
N.~Dunford and J.~T. Schwartz, {\em Spectral theory: self adjoint operators in
  Hilbert space}.
\newblock Interscience publishers, 1963.

\bibitem{lancaster1972norms}
P.~Lancaster and H.~K. Farahat, ``Norms on direct sums and tensor products,''
  {\em mathematics of computation}, vol.~26, no.~118, pp.~401--414, 1972.

\end{thebibliography}
\end{document}